\documentclass[runningheads]{llncs}
\usepackage[T1]{fontenc}
\usepackage{graphicx}
\usepackage{amsfonts,amssymb,amsmath}
\begin{document}
\title{Determinantal expressions of certain \\
integrals on symmetric spaces}
%
%
\author{Salem Said\inst{1} \and Cyrus Mostajeran\inst{2}
}
\authorrunning{Salem Said and Cyrus Mostajeran}
%
\institute{CNRS, Laboratoire Jean Kuntzmann (UMR 5224) \and
School of Physical and Mathematical Sciences, NTU Singapore}
%
\maketitle              
\begin{abstract}
The integral of a function $f$ defined on a symmetric space $M \simeq G/K$ may be expressed in the form of a determinant (or Pfaffian),\linebreak when $f$ is $K$-invariant and, in a certain sense, a tensor power of a positive function of a single variable. The paper presents a few examples of this idea and discusses future extensions. Specifically, the examples involve symmetric cones, Grassmann manifolds, and classical domains. 

\end{abstract}
%
%

\section{Introduction}
Riemannian symmetric spaces were classified by \'E. Cartan, back in the 1920s. A comprehensive account of this classification may be found in the monograph~\cite{helgasson}. In the 1960s, a classification of quantum symmetries led Dyson to introduce three kinds of random matrix ensembles, orthogonal, unitary, and symplectic~\cite{dyson}. These three kinds of ensembles are closely related to the symmetric spaces known as symmetric cones, and also to their compact duals, which provide for so-called circular ensembles. More recently, Dyson's classification of quantum symmetries has been extended to free fermionic systems. It turned out that this extended classification is in on-to-one correspondance with Cartan's old classification of symmetric spaces~\cite{zirnbauer}. This correspondance has motivated the notion that the relationship between random matrices and symmetric spaces extends well beyond symmetric cones, and is of a  general nature (for example~\cite{edelman} or~\cite{tierz,said3}).

The present submission has a modest objective. It is to show how the integral of a function $f$, defined on a symmetric space $M \simeq G/K$, can be expressed in the form of a determinant or Pfaffan, when $f$ is $K$-invariant and satisfies an additional hypothesis, formulated in Section \ref{sec:idea} below. This is not carried out in a general setting, but through a non-exhaustive set of examples, including symmetric cones, Grassmann manifolds, classical domains, and their duals (for the case of compact Lie groups, yet another example of symmetric spaces, see~\cite{meckes}).  

The determinantal expressions obtained here, although elementary, are an analytic pre-requisite to developing the \textit{random matrix theory of Riemannian symmetric spaces}. This long-term goal is the motivation behind the present work.
 
Unfortunately, due to limited space, no proofs are provided for statements made in the following. These will be given in an upcoming extended version.

\section{Integral formulas}
Let $M$ be a Riemannian symmetric space, given by the symmetric pair $(G,K)$. Write $\mathfrak{g} = \mathfrak{k} + \mathfrak{p}$ the corresponding Cartan decomposition, and let $\mathfrak{a}$ be a maximal abelian subspace of $\mathfrak{p}$. Then, denote by $\Delta$ a set of positive reduced roots on $\mathfrak{a}$~\cite{helgasson}. 

Assume that $\mathfrak{g} = \mathfrak{z}(\mathfrak{g}) + \mathfrak{g}_{\scriptscriptstyle ss}$ where $\mathfrak{z}(\mathfrak{g})$ is the centre of $\mathfrak{g}$ and $\mathfrak{g}_{\scriptscriptstyle ss}$ is semisimple and non-compact ($\mathfrak{g}_{\scriptscriptstyle ss}$ is a real Lie algebra). The Riemannian exponential $\mathrm{Exp}$ maps $\mathfrak{a}$ isometrically onto a totally flat submanifold of $M$, and any $x \in M$ is of the form
 $x = k\cdot \mathrm{Exp}(a)$ where $k \in K$ and $a \in \mathfrak{a}$.






Let $f:M \rightarrow \mathbb{R}$ be a $K$-invariant function, $f(k\cdot x) = f(x)$ for $k\in K$ and $x \in M$. There is no ambiguity in writing $f(x) = f(a)$ where $x = k\cdot \mathrm{Exp}(a)$. With this notation, there exists a constant $C_{\scriptscriptstyle M}$ such that~\cite{helgasson}
\begin{equation} \label{eq:noncompactintegral}
  \int_M f(x)\hspace{0.02cm}\mathrm{vol}(dx) = C_{\scriptscriptstyle M}\hspace{0.03cm}\int_{\mathfrak{a}}f(a)\prod_{\lambda\in\Delta} \sinh^{m_{\scriptscriptstyle \lambda}}\!\left|\lambda(a)\right|da
\end{equation}
where $da$ is the Lebesgue measure on $\mathfrak{a}$. 

 The dual $\hat{M}$ of $M$ is a symmetric space given by the symmetric pair $(U,K)$, where $U$ is a compact Lie group, with the Cartan decomposition $\mathfrak{u} = \mathfrak{k} + \mathrm{i}\hspace{0.02cm}\mathfrak{p}$\hfill\linebreak ($\mathrm{i} = \sqrt{-1}$). Now, $\mathrm{Exp}$ maps $\mathrm{i}\hspace{0.02cm}\mathfrak{a}$ onto a torus $T$ which is totally flat in $\hat{M}$, and any point $x \in \hat{M}$ is of the form $x = k\cdot \mathrm{Exp}(\mathrm{i}a)$ where $k \in K$ and $a \in \mathfrak{a}$.


If $f:\hat{M}\rightarrow \mathbb{R}$ is $K$-invariant, there is no ambiguity in writing $f(x) = f(t)$ where $x = k\cdot t$, $t = \mathrm{Exp}(\mathrm{i}a)$. In this notation~\cite{helgasson}, 
\begin{equation} \label{eq:compactintegral}
  \int_{\hat{M}} f(x)\hspace{0.02cm}\mathrm{vol}(dx) = C_{\scriptscriptstyle M}\hspace{0.03cm}\int_{T}f(t)\prod_{\lambda\in\Delta} \sin^{m_{\scriptscriptstyle \lambda}}\!\left|\lambda(t)\right|dt
\end{equation}
where $dt$ is the Haar measure on $T$. Here, $\sin\!\left|\lambda(t)\right| = \sin\!\left|\lambda(a)\right|$ where $t = \mathrm{Exp}(\mathrm{i}a)$, and this does not depend on the choice of $a$. 


\section{Determinantal expressions}
Let $\mu$ be a positive measure on a real interval $I$. Consider the multiple integrals, 
\begin{equation} \label{eq:zbetanc}
   z_\beta(\mu) = \frac{1}{N!}   \int_I\ldots\int_I\,  \left|V(u_{\scriptscriptstyle _1},\ldots,u_{\scriptscriptstyle N})\right|^\beta\,\mu(du_{\scriptscriptstyle 1})\ldots\mu(du_{\scriptscriptstyle N})
\end{equation}
where $V$ denotes the Vandermonde determinant and $\beta = 1,2$ or $4$. Consider also the following bilinear forms, 
\begin{equation} \label{eq:skew1}
 \left(h,g\hspace{0.02cm}\right)_{(\mu,1)} = \int_I\int_I (h(u)\hspace{0.02cm}\varepsilon(u-v)\hspace{0.02cm}g(v))\,\mu(du)\mu(dv)
\end{equation}
\begin{equation} \label{eq:scal2}
 \left(h,g\hspace{0.02cm}\right)_{(\mu,2)} = \int_I h(u)g(u)\,\mu(du)
\end{equation}
\begin{equation} \label{eq:skew4}
 \left(h,g\hspace{0.02cm}\right)_{(\mu,4)} = \int_I (h(u)g^\prime(u) - g(u)h^\prime(u))\,\mu(du)
\end{equation}
Here, $\varepsilon$ denotes the unit step function and the prime denotes the derivative. \hfill\linebreak
In the following proposition, $\det$ denotes the determinant and  $\mathrm{pf}$ the Pfaffian.
\begin{proposition} \label{prop:identities_rl}
The following hold for any probability measure $\mu$ as above. \\[0.1cm]
(a) if $N$ is even,
\begin{equation} \label{eq:debruijnone_even}
z_{\scriptscriptstyle 1}(\mu) = \mathrm{pf}\left\lbrace \left(u^k,u^\ell\hspace{0.02cm}\right)_{(\mu,1)}\right\rbrace^{\!N-1}_{\!k,\ell=0}
\end{equation}
(b) on the other hand, if $N$ is odd,
\begin{equation} \label{eq:debruijnone_odd}
  z_{\scriptscriptstyle 1}(\mu) = \mathrm{pf}\left\lbrace\! \begin{array}{rc} \left(u^k,u^\ell\hspace{0.02cm}\right)_{(\mu,1)} & \left(1,u^k\hspace{0.02cm}\right)_{(\mu,2)} \\[0.2cm] -\left(u^\ell\hspace{0.02cm},1\right)_{(\mu,2)} & 0\end{array}\!\right\rbrace^{\!N-1}_{\!k,\ell=0} 
\end{equation}
(c) moreover,
\begin{equation} \label{eq:andreev}
  z_{\scriptscriptstyle 2}(\mu) = \det  \left\lbrace \left(u^k,u^\ell\hspace{0.02cm}\right)_{(\mu,2)}\right\rbrace^{\!N-1}_{\!k,\ell=0}
\end{equation}
(d) and, finally,
\begin{equation} \label{eq:debruijnfour}
  z_{\scriptscriptstyle 4}(\mu) = \mathrm{pf}  \left\lbrace \left(u^k,u^\ell\hspace{0.02cm}\right)_{(\mu,4)}\right\rbrace^{\!2N-1}_{\!k,\ell=0}
\end{equation}
\vspace{0.1cm} 
\end{proposition}

On the other hand, if $\mu$ is a probability measure on the unit circle $S^1$, and 
\begin{equation} \label{eq:zbetac}
   z_\beta(\mu) = \frac{1}{N!}   \int_{S^1}\ldots\int_{S^1}\,  \left|V(u_{\scriptscriptstyle _1},\ldots,u_{\scriptscriptstyle N})\right|^\beta\,\mu(du_{\scriptscriptstyle 1})\ldots\mu(du_{\scriptscriptstyle N})
\end{equation}
consider the bilinear form
\begin{equation} \label{eq:skew1_circ}
 \left( h,g\hspace{0.02cm}\right)_{(\mu,1)} = \int^{2\pi}_0\int^{2\pi}_0 (h(e^{ix})\hspace{0.02cm}\varepsilon(x-y)\hspace{0.02cm}g(e^{iy}))\,\tilde{\mu}(dx)\tilde{\mu}(dy)
\end{equation}
where $\tilde{\mu}$ is the pullback of the measure $\mu$ through the map that takes $x$ to $e^{ix}$, and let $\left(h,g\hspace{0.02cm}\right)_{(\mu,2)}$ and $\left(h,g\hspace{0.02cm}\right)_{(\mu,4)}$ be given as in (\ref{eq:scal2}) and (\ref{eq:skew4}), with integrals over $S^1$ instead of $I$. 
\begin{proposition} \label{prop:identities_uc}
The following hold for any probability measure $\mu$ on $S^1$. \\[0.1cm]
(a) if $N$ is even,
\begin{equation} \label{eq:debruijnone_even_circ}
z_{\scriptscriptstyle 1}(\mu) = (-i)^{N(N-1)/2}\,\times\,\mathrm{pf}\left\lbrace \left(g_k,g_\ell\hspace{0.02cm}\right)_{(\mu,1)}\right\rbrace^{\!N-1}_{\!k,\ell=0}
\end{equation}
where $g_k(u) = u^{k-(N-1)/2}\,$. \\[0.15cm]
(b) on the other hand, if $N$ is odd,
\begin{equation} \label{eq:debruijnone_odd_circ}
  z_{\scriptscriptstyle 1}(\mu) = (-i)^{N(N-1)/2}\,\times\,\mathrm{pf}\left\lbrace\! \begin{array}{rc} \left(g_k,g_\ell\hspace{0.02cm}\right)_{(\mu,1)} & \left(1,g_k\hspace{0.02cm}\right)_{(\mu,2)} \\[0.2cm] -\left(g_\ell\hspace{0.02cm},1\right)_{(\mu,2)} & 0\end{array}\!\right\rbrace^{\!N-1}_{\!k,\ell=0} 
\end{equation}
with the same definition of $g_k(u)$. \\[0.15cm]
(c) moreover,
\begin{equation} \label{eq:andreev_circ}
  z_{\scriptscriptstyle 2}(\mu) = \det  \left\lbrace \left(u^k,u^{-\ell}\hspace{0.02cm}\right)_{(\mu,2)}\right\rbrace^{\!N-1}_{\!k,\ell=0}
\end{equation}
(d) and, finally,
\begin{equation} \label{eq:debruijnfour_circ}
  z_{\scriptscriptstyle 4}(\mu) = \mathrm{pf}  \left\lbrace \left(h_k,h_\ell\hspace{0.02cm}\right)_{(\mu,4)}\right\rbrace^{\!2N-1}_{\!k,\ell=0}
\end{equation}
\vspace{0.1cm} 
where $h_k(u) = u^{k-(N-1)}$. 
\end{proposition}
Both of the above Propositions \ref{prop:identities_rl} and \ref{prop:identities_uc} are directly based on~\cite{mehta}. 

\section{Main idea} \label{sec:idea}
An additional hypothesis is made on the function $f(a)$ (in (\ref{eq:noncompactintegral})) or $f(t)$ (in (\ref{eq:compactintegral}))\,: that there exists a natural orthonormal basis $(e_j;j=1,\ldots,r)$ of $\mathfrak{a}$, such that 
\begin{equation} \label{eq:hypothesis}
f(a) = \prod^r_{j=1}w(a_j) \hspace{1cm} f(t) = \prod^r_{j=1} w(t_j)
\end{equation}
where $w$ is a positive function of a single variable, and $a_j$ are the components of $a$ in the basis $(e_j;j=1,\ldots,r)$, while $t_j = \mathrm{Exp}(\mathrm{i}\hspace{0.02cm}a_je_j)$. In this sense, it may be said that $f$ is the $r$-th tensor power of $w$. 

What is meant by \textit{natural} is that (\ref{eq:hypothesis}) will imply that the integral (\ref{eq:noncompactintegral}) or (\ref{eq:compactintegral}) can be transformed into a multiple integral of the form (\ref{eq:zbetanc}) or (\ref{eq:zbetac}), respectively. Thus, in the case of (\ref{eq:noncompactintegral}), there exists a measure $\mu$ on an interval $I$, which satisfies
$$
  \int_M f(x)\hspace{0.02cm}\mathrm{vol}(dx) = \tilde{C}_{\scriptscriptstyle M}\times z_\beta(\mu) \hspace{1cm} \text{($\tilde{C}_{\scriptscriptstyle M}$ is a new constant)}
$$
and, in the case of (\ref{eq:compactintegral}), there is a measure $\mu$ on $S^1$, which yields a similar identity. It should be noted that this measure $\mu$ will depend on the function $w$ from (\ref{eq:hypothesis}).

Then, Propositions \ref{prop:identities_rl} and \ref{prop:identities_uc} provide a determinantal (or Pfaffian) expression of the initial integral on the symmetric space $M$ or $\hat{M}$.

At present, this is not a theorem, but a mere idea or observation, supported by the examples in the following section.

\section{Examples}
\subsection{Symmetric cones} \label{subsec:cones}
Consider the following Lie groups (in the usual notation, as found in~\cite{helgasson}).
$$
\begin{array}{lllllll}
\beta & \hspace{0.4cm} & G_\beta &\hspace{0.2cm}& U_\beta &\hspace{0.2cm}& K_\beta \\[0.2cm]
1      & & GL_N(\mathbb{R}) & &U(N) && O(N) \\[0.15cm]
2      & & GL_N(\mathbb{C}) & & U(N) \times U(N) && U(N) \\[0.15cm]
4      & & GL_N(\mathbb{H}) & & U(2N) && Sp(N)
\end{array}
$$
Then, $M_\beta \simeq G_\beta/K_\beta$ is a Riemannian symmetric space, with dual $\hat{M}_\beta = U_\beta/K_\beta\hspace{0.03cm}$. In fact, $M_\beta$ is realised as a so-called symmetric cone\,: the cone of positive-definite real, complex, or quaternion matrices (according to the value of $\beta=1,2$ or $4$). 

Each $x \in M_\beta$ is of the form $k\lambda\,k^\dagger$ where $k \in K_\beta$ and $\lambda$ is a positive diagonal matrix ($\dagger$ denotes the transpose, conjugate-transpose, or quaternion conjugate-transpose). If $f:M_\beta \rightarrow \mathbb{R}$ is $K_\beta$-invariant, and can be written $f(x) = \prod w(\lambda_j)$, 
\begin{equation} \label{eq:example_cones}
   \int_{M_\beta} f(x)\hspace{0.02cm}\mathrm{vol}(dx) = \tilde{C}_{\beta}\times z_\beta(\mu)
\end{equation}
where $\mu(du) = (w(u)u^{-\scriptscriptstyle N_\beta})\hspace{0.02cm}du$, with $N_\beta = (\beta/2)(N-1)+1$, on the interval $I = (0,\infty)$. The constant $\tilde{C}_\beta$ is known explicitly, but this is irrelevant at present. 

The dual $\hat{M}_\beta$ can be realised as the space of symmetric unitary matrices ($\beta = 1$), of unitary matrices ($\beta = 2$), or of antisymmetric unitary matrices with double dimension $2N$, ($\beta = 4$). 

If $\beta = 1,2$, then $x \in \hat{M}_\beta$ is of the form $ke^{\mathrm{i}\theta}\,k^\dagger$ where $k \in K_\beta$ and $\theta$ is real diagonal. However, if $\beta = 4$, there is a somewhat different matrix factorisation,
\begin{equation} \label{eq:symplecdiag}
  x = k\,\left(\begin{array}{ll} & -e^{\mathrm{i}\theta} \\[0.1cm]
e^{\mathrm{i}\theta} & \end{array}\right) k^{\mathrm{tr}} \hspace{1cm} \text{($\mathrm{tr}$ denotes the transpose)}
\end{equation}
where $k \in Sp(N)$ is considered as a $2N\times 2N$ complex matrix (rather than a $N\times N$ quaternion matrix). If $f:\hat{M}_\beta \rightarrow \mathbb{R}$ is $K_\beta$-invariant, $f(x) = \prod w(e^{\mathrm{i}\theta_j})$,\hfill\linebreak
\begin{equation} \label{eq:example_unitaries}
   \int_{\hat{M}_\beta} f(x)\hspace{0.02cm}\mathrm{vol}(dx) = \tilde{C}_{\beta}\times z_\beta(\mu)
\end{equation}
where $\mu(du) = w(u)|du|$ on the unit circle $S^1$ ($|du| = d\varphi$ if $u = e^{\mathrm{i}\varphi}$).\\[0.1cm]
\noindent \textbf{Remark\,:} in many textbooks, $\hat{M}_{\scriptscriptstyle 1}$ is realised as the space of real structures on $\mathbb{C}^N$, and $\hat{M}_{\scriptscriptstyle 4}$ as the space of quaternion structures on $\mathrm{C}^{2N}$. The alternative realisations proposed here seem less well-known, but more concrete, so to speak.

\subsection{Grassmann manifolds}
Consider the following Lie groups (again, for the notation, see~\cite{helgasson}).
$$
\begin{array}{lllllll}
\beta & \hspace{0.4cm} & G_\beta &\hspace{0.2cm}& U_\beta &\hspace{0.2cm}& K_\beta \\[0.2cm]
1      & & O(p,q) & & O(p+q) && O(p)\times O(q) \\[0.15cm]
2      & & U(p,q) & & U(p+q) && U(p)\times U(q) \\[0.15cm]
4      & & Sp(p,q) & & Sp(p+q) && Sp(p) \times Sp(q)
\end{array}
$$
Then, $M_\beta \simeq G_\beta/K_\beta$ is a Riemannian symmetric space, with dual $\hat{M}_\beta = U_\beta/K_\beta\hspace{0.03cm}$. 
The $M_\beta$ may be realised as follows~\cite{huang} ($\mathbb{K} = \mathbb{R},\mathbb{C}$ or $\mathbb{H}$, according to $\beta$),
\begin{equation} \label{eq:spacelike}
  M_\beta = \lbrace x : \text{$x$ is a $p$-dimensional and space-like subspace of $\mathbb{K}^{p+q}\rbrace$}
\end{equation}
Here, $x$ is space-like if $|\xi_p|^2-|\xi_q|^2 > 0$ for all $\xi \in x$ with $\xi = (\xi_p\hspace{0.02cm},\xi_q)$, where $|\cdot|$ denotes the standard Euclidean norm on $\mathbb{K}^p$ or $\mathbb{K}^q$. Moreover, for each $x \in M_{\beta}\hspace{0.02cm}$, $x = k(x_\tau)$ where $k \in K_\beta$ and $x_\tau \in M_\beta$ is spanned by the vectors 
$$
\cosh(\tau_j)\xi_j + \sinh(\tau_j)\xi_{p+j} \hspace{1cm} j = 1,\ldots, p
$$ 
with $(\xi_k;k=1,\ldots,p+q)$  the canonical basis of $\mathbb{K}^{p+q}$, and $(\tau_j;j=1,\ldots, p)$ real ($p \leq q$ throughout this paragraph).

If $f:M_\beta\rightarrow \mathbb{R}$ is $K_\beta$-invariant, $f(x) =f(\tau)$, the right-hand side of (\ref{eq:noncompactintegral}) reads
(the positive reduced roots can be found in~\cite{sakai})
\begin{equation} \label{eq:spacelikeintegral}
C_{\beta}\hspace{0.03cm}\int_{\mathbb{R}^p}f(\tau)\hspace{0.02cm}
           \prod^p_{j=1}\sinh^{\beta(q-p)}|\tau_j|\hspace{0.02cm}\sinh^{\beta-1}|2\tau_j|\prod_{i<j}\left|\cosh(2\tau_i) - \cosh(2\tau_j)\right|^\beta\hspace{0.02cm}
d\tau
\end{equation}
and this can be transformed into the form (\ref{eq:zbetanc}), by introducing $u_j = \cosh(2\tau_j)$. This will reappear, with $\beta = 2$ and $p = q$, in the following paragraph. 

Now, the duals $\hat{M}_\beta$ are real, complex, or quaternion Grassmann manifolds, 
\begin{equation} \label{eq:grassmann}
  \hat{M}_\beta = \lbrace x : \text{$x$ is a $p$-dimensional subspace of $\mathbb{K}^{p+q}\rbrace$}
\end{equation}
For each $x \in \hat{M}_\beta\hspace{0.02cm}$, $x = k(x_\theta)$ where $k \in K_\beta$ and $x_\theta$ is spanned by the vectors 
$$
\cos(\theta_j)\xi_j + \sin(\theta_j)\xi_{p+j} \hspace{1cm} j = 1,\ldots, p
$$
with $(\theta_j;j=1,\ldots, p)$ real. 

If $f:\hat{M}_\beta\rightarrow \mathbb{R}$ is $K_\beta$-inariant, $f(x) =f(\theta)$, the right-hand side of (\ref{eq:compactintegral}) reads
\begin{equation} \label{eq:grassmannintegral}
C_{\beta}\hspace{0.03cm}\int_{(0,\pi)^p}f(\theta)\hspace{0.02cm}
           \prod^p_{j=1}\sin^{\beta(q-p)}|\theta_j|\hspace{0.02cm}\sin^{\beta-1}|2\theta_j|\prod_{i<j}\left|\cos(2\theta_i) - \cos(2\theta_j)\right|^\beta\hspace{0.02cm}
d\theta
\end{equation}
which can be transformed into the form (\ref{eq:zbetac}), by introducing $u_j = \cos(2\theta_j)$. In~\cite{edelman}, this is used to recover the Jacobi ensembles of random matrix theory.  \\[0.1cm]
\noindent \textbf{Remark\,:} the angles $\theta_j$ may be taken in the interval $(-\pi/2,\pi/2)$ instead of $(0,\pi)$. In this case, $|\theta_j|$ are the principal angles between $x_\theta$ and the subspace $x_o$ spanned by $(\xi_j;j=1,\ldots,p)$. By analogy, it is natural to think of $|\tau_j|$ as the `principal boosts' (using the language of special relativity) between $x_\tau$ and $x_o\hspace{0.02cm}$.

\subsection{Classical domains}
Consider, finally, the following Lie groups (again, for the notation, see~\cite{helgasson}).
$$
\begin{array}{lllllll}
\beta & \hspace{0.4cm} & G_\beta &\hspace{0.2cm}& U_\beta &\hspace{0.2cm}& K_\beta \\[0.2cm]
1      & & Sp(N,\mathbb{R}) & & Sp(N) && U(N) \\[0.15cm]
2      & & U(N,N) & & U(2N) && U(N)\times U(N) \\[0.15cm]
4      & & O^*(4N) & & O(4N) && U(2N) 
\end{array}
$$
Then, $M_\beta \simeq G_\beta/K_\beta$ is a Riemannian symmetric space, with dual $\hat{M}_\beta = U_\beta/K_\beta\hspace{0.03cm}$. The $M_\beta$ are realised as classical domains, whose elements are $N\times N$ complex matrices (if $\beta = 1, 2$) or $2N \times 2N$ complex matrices (if $\beta = 4$), with operator norm $< 1$, and which are in addition symmetric ($\beta = 1$) or antisymmetric ($\beta =4$). 

If $\beta = 1,2$, then any $x \in M_\beta$ may be written
\begin{equation} \label{eq:takagi}
x = k_{\scriptscriptstyle 1}\hspace{0.02cm}(\tanh(\lambda))\hspace{0.02cm}k_{\scriptscriptstyle 2}
\end{equation}
where $k_{\scriptscriptstyle 1}$ and $k_{\scriptscriptstyle 2}$ are unitary ($k^{\phantom{\mathrm{tr}}}_{\scriptscriptstyle 2} = k^\mathrm{tr}_{\scriptscriptstyle 1}$, in case $\beta = 1$), and $\lambda$ is real diagonal. However, if $\beta = 4$,
\begin{equation} \label{eq:skewsiegel}
x =  k\,\left(\begin{array}{ll} & -\tanh(\lambda) \\[0.1cm]
\tanh(\lambda) & \end{array}\right) k^{\mathrm{tr}}
\end{equation}
where $k$ is $2N \times 2N$ unitary. If $f:M_\beta \rightarrow \mathbb{R}$ is $K_\beta$-invariant, and $f(x) = \prod w(\lambda_j)$,
$$
    \int_{M_\beta} f(x)\mathrm{vol}(dx) = \tilde{C}_{\beta}\hspace{0.02cm}\int_{\mathbb{R}^N}\prod^N_{j=1} w(\lambda_j)\sinh|2\lambda_j|\prod_{i<j}\left|\cosh(2\lambda_i) - \cosh(2\lambda_j)\right|^\beta
d\lambda
$$
After introducing $u_j = \cosh(2\lambda_j)$, this immediately becomes
\begin{equation} \label{eq:example_domains}
   \int_{M_\beta} f(x)\mathrm{vol}(dx) =  \tilde{C}_{\beta}\times z_\beta(\mu)
\end{equation}
where $\mu(du) = w(\mathrm{acosh}(u)/2)\hspace{0.02cm}du$ on the interval $I = (1,\infty)$. \\[0.1cm]
\textbf{Remark\,:} the domain $M_{2}$ is sometimes called the Siegel disk. As an application of (\ref{eq:example_domains}), consider a random $x \in M_{2}$ with a Gaussian probability density function
\begin{equation}    \label{eq:siegelgaus}
   p(x|\bar{x},\sigma) = \left(Z(\sigma)\right)^{-1}\hspace{0.02cm}\exp\left[-\frac{d^{\hspace{0.02cm} 2}(x\hspace{0.02cm},\bar{x})}{2\sigma^2}\right]
\end{equation}
with respect to $\mathrm{vol}(dx)$, where $d(x\hspace{0.02cm},\bar{x})$ denotes Riemannian distance and $\sigma > 0$. Then, following the arguments in~\cite{said3}, (\ref{eq:example_domains}) can be used to obtain
$$
  Z(\sigma) = \tilde{C}_{\scriptscriptstyle 2}\times \det \left\lbrace m_{k+\ell}(\sigma)\right\rbrace^{N-1}_{k,\ell=0} \hspace{0.25cm} m_{j}(\sigma) = \int^{\infty}_{1}\hspace{0.02cm} \exp\left(-\mathrm{acosh}^2(u)/8\sigma^2\right)u^j\hspace{0.03cm}du
$$
The integrals $m_j(\sigma)$ are quite easy to compute, and one is then left with a determinantal expression of $Z(\sigma)$. The starting point to the study of the random matrix $x$ is the following observation. If $x$ is written as in (\ref{eq:takagi}) and $u_j = \cosh(2\lambda_j)$, then the random subset $\lbrace u_j;j=1,\ldots,N\rbrace$ of $I=(1,\infty)$ is a determinantal point process (see~\cite{kj}). By writing down its kernel function, one may begin to investigate in detail many of its statistical properties, including asymptotic ones, such as the asymptotic density of the $(u_j)$, or the asymptotic distribution of their maximum, in the limit where $N \rightarrow \infty$ (of course, with suitable re-scaling).
\vfill
\pagebreak
\section{Future directions}
The present submission developed determinantal expressions for integrals on symmetric spaces on a case-by-case basis, only through a non-exhaustive set of examples. Future work should develop these expressions in a fully general way, by transforming (\ref{eq:noncompactintegral}) and (\ref{eq:compactintegral}) into (\ref{eq:zbetanc}) or (\ref{eq:zbetac}), for any system of reduced roots.

The long-term goal is to understand the \textit{random matrix theory of symmetric spaces}. One aspect of this is to understand the asymptotic properties of a joint probability density (in the notation of (\ref{eq:noncompactintegral}))
$$
f(a)\prod_{\lambda\in\Delta} \sinh^{m_{\scriptscriptstyle \lambda}}\!\left|\lambda(a)\right|da
$$
and analyse how these depend on the set of positive reduced roots $\Delta$. It is worth mentioning that, in previous work~\cite{said3}, it was seen that a kind of universality holds, where different root systems lead to the same asymptotic properties.

Random matrix theory (in its classical realm of orthogonal, unitary, and symplectic ensembles) has so many connections to physics, combinatorics, and complex systems in general. A further important direction is to develop such connections for the random matrix theory of symmetric spaces.

\end{document}